\documentclass{amsart}
\usepackage{enumerate}
\usepackage{amsmath,amssymb,amsfonts,amscd,amsthm,amsopn,epsfig}  
\usepackage{amsthm}
\usepackage{amssymb}
\usepackage{latexsym}
\usepackage{url}

\usepackage{hyperref}
\hypersetup{
    citecolor=blue,        
}

\usepackage{color}

\newtheorem{definition}{Definition}[section]

\newtheorem{rmk}{Remark}[section]

\numberwithin{equation}{section}

\numberwithin{equation}{section}
\def\cA{{\cal A}}   \def\cB{{\cal B}}   
\def\cD{{\cal D}}     
      
      \def\cL{{\cal L}}
      \def\cO{{\cal O}}
      \def\cR{{\cal R}}
   \def\cT{{\cal T}}   \def\cU{{\cal U}}
      
\def\cY{{\cal Y}}

\newcommand{\beqa}{\begin{eqnarray}}
\newcommand{\eeqa}{\end{eqnarray}}

\newcommand{\nc}{\newcommand}
\newcommand{\rnc}{\renewcommand}

\nc{\cal}{\mathcal}

\nc{\goth}{\mathfrak}
\rnc{\bold}{\mathbf}
\renewcommand{\frak}{\mathfrak}
\renewcommand{\Bbb}{\mathbb}

\nc{\Cal}{\mathcal}
\nc{\Xp}[1]{X^+(#1)}
\nc{\Xm}[1]{X^-(#1)}
\nc{\on}{\operatorname}
\nc{\ch}{\mbox{ch}}
\nc{\Z}{{\bold Z}}
\nc{\J}{{\mathcal J}}
\nc{\C}{{\bold C}}
\nc{\Q}{{\bold Q}}
\nc{\oC}{{\widetilde{C}}}
\nc{\oc}{{\tilde{c}}}
\nc{\og}{{\tilde{\gamma}}}
\nc{\lC}{{\overline{C}}}
\nc{\lc}{{\overline{c}}}
\nc{\Rt}{{\tilde{R}}}

\nc{\odel}{{\overline{\delta}}}

\nc{\N}{{\Bbb N}}
\nc\beq{\begin{equation}}
\nc\enq{\end{equation}}
\nc\lan{\langle}
\nc\ran{\rangle}
\nc\bsl{\backslash}
\nc\mto{\mapsto}
\nc\lra{\leftrightarrow}
\nc\hra{\hookrightarrow}
\nc\sm{\smallmatrix}
\nc\esm{\endsmallmatrix}
\nc\sub{\subset}
\nc\ti{\tilde}
\nc\nl{\newline}
\nc\fra{\frac}
\nc\und{\underline}
\nc\ov{\overline}
\nc\ot{\otimes}
\nc\bbq{\bar{\bq}_l}
\nc\bcc{\thickfracwithdelims[]\thickness0}
\nc\ad{\text{\rm ad}}
\nc\Ad{\text{\rm Ad}}
\nc\Hom{\text{\rm Hom}}
\nc\End{\text{\rm End}}
\nc\Ind{\text{\rm Ind}}
\nc\Res{\text{\rm Res}}
\nc\Ker{\text{\rm Ker}}
\rnc\Im{\text{Im}}
\nc\sgn{\text{\rm sgn}}
\nc\tr{\text{\rm tr}}
\nc\Tr{\text{\rm Tr}}
\nc\supp{\text{\rm supp}}
\nc\card{\text{\rm card}}
\nc\bst{{}^\bigstar\!}
\nc\he{\heartsuit}
\nc\clu{\clubsuit}
\nc\spa{\spadesuit}
\nc\di{\diamond}
\nc\al{\alpha}
\nc\bet{\beta}
\nc\ga{\gamma}
\nc\de{\delta}
\nc\ep{\epsilon}
\nc\io{\iota}
\nc\om{\omega}
\nc\si{\sigma}
\rnc\th{\theta}
\nc\ka{\kappa}
\nc\la{\lambda}
\nc\ze{\zeta}

\nc\vp{\varpi}
\nc\vt{\vartheta}
\nc\vr{\varrho}

\nc\Ga{\Gamma}
\nc\De{\Delta}
\nc\Om{\Omega}
\nc\Si{\Sigma}
\nc\Th{\Theta}
\nc\La{\Lambda}

\nc\boa{\bold a}
\nc\bob{\bold b}
\nc\boc{\bold c}
\nc\bod{\bold d}
\nc\boe{\bold e}
\nc\bof{\bold f}
\nc\bog{\bold g}
\nc\boh{\bold h}
\nc\boi{\bold i}
\nc\boj{\bold j}
\nc\bok{\bold k}
\nc\bol{\bold l}
\nc\bom{\bold m}
\nc\bon{\bold n}
\nc\boo{\bold o}
\nc\bop{\bold p}
\nc\boq{\bold q}
\nc\bor{\bold r}
\nc\bos{\bold s}
\nc\bou{\bold u}
\nc\bov{\bold v}
\nc\bow{\bold w}
\nc\boz{\bold z}

\nc\ba{\bold A}
\nc\bb{\bold B}
\nc\bc{\bold C}
\nc\bd{\bold D}
\nc\be{\bold E}
\nc\bg{\bold G}
\nc\bh{\bold H}
\nc\bi{\bold I}
\nc\bj{\bold J}
\nc\bl{\bold L}
\nc\bm{\bold M}
\nc\bn{\bold N}
\nc\bo{\bold O}
\nc\bp{\bold P}
\nc\bq{\bold Q}
\nc\br{\bold R}
\nc\bs{\bold S}
\nc\bt{\bold T}
\nc\bu{\bold U}
\nc\bv{\bold V}
\nc\bw{\bold W}
\nc\bz{\bold Z}
\nc\bx{\bold X}

\nc\ca{\mathcal A}
\nc\cb{\mathcal B}
\nc\cc{\mathcal C}
\nc\cd{\mathcal D}
\nc\ce{\mathcal E}
\nc\cf{\mathcal F}
\nc\cg{\mathcal G}
\rnc\ch{\mathcal H}
\nc\ci{\mathcal I}
\nc\cj{\mathcal J}
\nc\ck{\mathcal K}
\nc\cl{\mathcal L}
\nc\cm{\mathcal M}
\nc\cn{\mathcal N}
\nc\co{\mathcal O}
\nc\cp{\mathcal P}
\nc\cq{\mathcal Q}
\nc\car{\mathcal R}
\nc\cs{\mathcal S}
\nc\ct{\mathcal T}
\nc\cu{\mathcal U}
\nc\cv{\mathcal V}
\nc\cz{\mathcal Z}
\nc\cx{\mathcal X}
\nc\cy{\mathcal Y}

\nc\bH{\textbf{H}}
\nc\bX{\textbf{X}}
\nc\bA{\textbf{A}}
\nc\bB{\textbf{B}}
\nc\bK{\textbf{K}}
\nc\bJ{\textbf{J}}
\nc\bE{\textbf{E}}
\nc\bF{\textbf{F}}
\nc\bk{\textbf{k}}

\nc\e[1]{E_{#1}}
\nc\ei[1]{E_{\delta - \alpha_{#1}}}
\nc\esi[1]{E_{s \delta - \alpha_{#1}}}
\nc\eri[1]{E_{r \delta - \alpha_{#1}}}
\nc\ed[2][]{E_{#1 \delta,#2}}
\nc\ekd[1]{E_{k \delta,#1}}
\nc\emd[1]{E_{m \delta,#1}}
\nc\erd[1]{E_{r \delta,#1}}

\nc\ef[1]{F_{#1}}
\nc\efi[1]{F_{\delta - \alpha_{#1}}}
\nc\efsi[1]{F_{s \delta - \alpha_{#1}}}
\nc\efri[1]{F_{r \delta - \alpha_{#1}}}
\nc\efd[2][]{F_{#1 \delta,#2}}
\nc\efkd[1]{F_{k \delta,#1}}
\nc\efmd[1]{F_{m \delta,#1}}
\nc\efrd[1]{F_{r \delta,#1}}

\nc\fa{\frak a}
\nc\fb{\frak b}
\nc\fc{\frak c}
\nc\fd{\frak d}
\nc\fe{\frak e}
\nc\ff{\frak f}
\nc\fg{\frak g}
\nc\fh{\frak h}
\nc\fj{\frak j}
\nc\fk{\frak k}
\nc\fl{\frak l}
\nc\fm{\frak m}
\nc\fn{\frak n}
\nc\fo{\frak o}
\nc\fp{\frak p}
\nc\fq{\frak q}
\nc\fr{\frak r}
\nc\fs{\frak s}
\nc\ft{\frak t}
\nc\fu{\frak u}
\nc\fv{\frak v}
\nc\fz{\frak z}
\nc\fx{\frak x}
\nc\fy{\frak y}

\nc\fA{\frak A}
\nc\fB{\frak B}
\nc\fC{\frak C}
\nc\fD{\frak D}
\nc\fE{\frak E}
\nc\fF{\frak F}
\nc\fG{\frak G}
\nc\fH{\frak H}
\nc\fJ{\frak J}
\nc\fK{\frak K}
\nc\fL{\frak L}
\nc\fM{\frak M}
\nc\fN{\frak N}
\nc\fO{\frak O}
\nc\fP{\frak P}
\nc\fQ{\frak Q}
\nc\fR{\frak R}
\nc\fS{\frak S}
\nc\fT{\frak T}
\nc\fU{\frak U}
\nc\fV{\frak V}
\nc\fZ{\frak Z}
\nc\fX{\frak X}
\nc\fY{\frak Y}
\nc\tfi{\ti{\Phi}}

\rnc\bol{\bold 1}

\newcommand{\hs}[1]{\hspace{#1 mm}}
\newcommand{\mb}[1]{\hs{4}\mbox{#1}\hs{4}}

\nc\ua{\bold U_\A}

\nc\RR{\mathbb R}
\nc\CC{\mathbb C}
\nc\II{\mathbb I}
\nc\qinti[1]{[#1]_i}
\nc\q[1]{[#1]_q}
\nc\xpm[2]{E_{#2 \delta \pm \alpha_#1}}  
\nc\xmp[2]{E_{#2 \delta \mp \alpha_#1}}
\nc\xp[2]{E_{#2 \delta + \alpha_{#1}}}
\nc\xm[2]{E_{#2 \delta - \alpha_{#1}}}
\nc\hik{\ed{k}{i}}
\nc\hjl{\ed{l}{j}}
\nc\qcoeff[3]{\left[ \begin{smallmatrix} {#1}& \\ {#2}& \end{smallmatrix}
\negthickspace \right]_{#3}}
\nc\qi{q}
\nc\qj{q}

\nc\ufdm{{_\ca\bu}_{\rm fd}^{\le 0}}

\nc\nonu{\nonumber\\}

\nc\isom{\cong} 

\nc{\pone}{{\Bbb C}{\Bbb P}^1}
\nc{\pa}{\partial}

\nc{\F}{{\mathcal F}}
\nc{\Sym}{{\goth S}}
\nc{\A}{{\mathcal A}}
\nc{\arr}{\rightarrow}
\nc{\larr}{\longrightarrow}

\nc{\ri}{\rangle}
\nc{\lef}{\langle}
\nc{\W}{{\mathcal W}}
\nc{\uqatwoatone}{{U_{q,1}}(\su)}
\nc{\uqtwo}{U_q(\goth{sl}_2)}
\nc{\dij}{\delta_{ij}}
\nc{\divei}{E_{\alpha_i}^{(n)}}
\nc{\divfi}{F_{\alpha_i}^{(n)}}
\nc{\Lzero}{\Lambda_0}
\nc{\Lone}{\Lambda_1}
\nc{\ve}{\varepsilon}
\nc{\phioneminusi}{\Phi^{(1-i,i)}}
\nc{\phioneminusistar}{\Phi^{* (1-i,i)}}
\nc{\phii}{\Phi^{(i,1-i)}}
\nc{\Li}{\Lambda_i}
\nc{\Loneminusi}{\Lambda_{1-i}}
\nc{\vtimesz}{v_\ve \otimes z^m}

\nc{\asltwo}{\widehat{\goth{sl}_2}}
\nc\ag{\widehat{\goth{g}}}  
\nc\teb{\tilde E_\boc}
\nc\tebp{\tilde E_{\boc'}}

\newcommand{\ben}{\begin{eqnarray}}
\newcommand{\een}{\end{eqnarray}}

\begin{document}

\title{Coideal algebras from twisted Manin triples}
\author{S. Belliard} 
\author{N. Cramp\'e}
\email{samuel.belliard@univ-montp2.fr, nicolas.crampe@univ-montp2.fr}
\address{Universit{\'e} Montpellier 2, Laboratoire Charles Coulomb UMR 5221, F-34095
Montpellier, France.
CNRS, Laboratoire Charles Coulomb, UMR 5221, F-34095 Montpellier.}
\begin{abstract}
We propose a new approach to study coideal algebras.
It is well-known that Manin triples (or equivalently Lie bi-algebra structures)
are the requirement to deform Lie algebras and to obtain quantum groups.
In this paper, introducing some particular automorphisms of Manin triples, we   
define new structures that we call Lie bi-ideal structures. A link with coisotropic subalgebras is explained. 
We show that their deformation 
provide coideal algebras. 
As examples, we recover from our general construction the twisted Yangians $\cY^{\pm}(sl(2))$, 
the q-Onsager algebra 
and the augmented q-Onsager algebra. As an important by-product, we find a new presentation for 
the twisted Yangians.
\end{abstract}

\maketitle

MSC numbers: 81R50, 17B37

Key words: Coideal; Coisotropic subalgebra; Deformation; Manin triple; Onsager algebras; Twisted Yangians;
Reflection algebras. 

\section{Introduction}

The quantum groups introduced in \cite{J1,Dr} have a lot of influence in 
different scientific domains, such as the study of the Hopf algebras 
or the development of the quantum inverse scattering method \cite{FST} to solve the quantum integrable systems. 
It was realized later that a Hopf structure is necessary to study the bulk part of the integrable systems
whereas a coideal structure is needed to study the boundaries.
Indeed, the study of integrable models with open boundary conditions (instead of periodic boundary conditions) 
requires the reflection equations \cite{Cher} which define subalgebras of quantum groups, 
called reflection algebras \cite{Sk} (see also \cite{MR}). These subalgebras are not
Hopf algebras anymore but are still coideal subalgebras. 

A lot of coideal subalgebras have been introduced:
the twisted quantized enveloping algebras corresponding to symmetric pairs $(gl_n,o_n)$ or 
$(gl_{2n},sp_{2n})$ \cite{No,NS,GK} or to symmetric pairs $(gl_{n+m},gl_n\oplus gl_m)$ \cite{DN,DNS},
the twisted Yangians  \cite{Ols,MNO}, 
the twisted q-Yangians \cite{MRS} or the q-Onsager algebras and their generalization \cite{IT1,IT2,BS,BB}
(see also \cite{DMS,DM} for an application to integrable quantum field theories). 
Let us also mention the papers \cite{Le} where a classification of coideals for the finite 
quantum groups has been performed, \cite{DKM, Mu} where Drinfeld twists 
allow them to introduce a universal reflection equation and \cite{KS} where they exhibit relations 
between the different approaches.

In this paper, we propose to modify the historical methods,
used by Drinfeld to introduce quantum groups \cite{Dr}, in order to obtain coideal subalgebras. 
Namely, Drinfeld used the Lie bi-algebra structure to characterize the expansion of the coproduct. 
Then, he showed that 
it defined completely the deformation and obtained quantum groups. 
In the case of coideal subalgebras, we have to replace
the notion of bi-algebra by a new structure: the Lie bi-ideal structure 
(see definition \ref{def:tau}). We also link this structure with the coisotropic subalgebra \cite{Lu,Ci}.
Then, we show how to deform this structure to get coideal algebras.  

One of the practical methods to get examples of Lie bi-algebra structures consists in finding Manin triples.
Indeed, there is a one-to-one correspondence between a Manin triple and a Lie bi-algebra structure (for finite 
algebras). In this paper, we associate a Lie bi-ideal structure to 
a particular automorphism of Manin triple, called the Manin triple twist. {}From a simple criterion on this twist, 
we identify two types of deformation, leading to Hopf subalgebras or coideal subalgebras.

The plan of this paper is as follows. In section \ref{sec:MT}, we recall the definition
of a Manin triple and its relation with the Lie bi-algebra structure. 
Then, we define the Manin triple twists and the Lie bi-ideal structure. 
In section \ref{sec:defo}, we are interested in the 
deformation of the previous structure. We recall
how to deform the Lie bi-algebra structure and explain how to deform 
the Lie bi-ideal structure.
In section \ref{sec:exM}, we give different examples of these structures associated to  
the loop algebra (section \ref{sec:lie}) and to the half-loop algebra (section \ref{sec:HLo}).
For each case, we provide a Lie bi-algebra structure and give the Hopf deformation to
get the quantum algebra and the Yangian.
We give also different examples of Manin triple twists and get the q-Onsager algebra,
the augmented q-Onsager algebra, the positive Borel of 
${\cU_\hbar(A_1^{(2)})}$ and the twisted Yangians ${\cY^\pm_\hbar(sl_2)}$.
In section \ref{sec:ReflEq}, we construct, using intertwiner techniques, scalar $K$-matrices 
from the coideal algebras obtained in the examples. It allows us to discuss the connection 
between our examples and some known reflection algebras. 
Finally, we conclude and give some open problems in section \ref{sec:conc}.

\section{Manin triple and twists \label{sec:MT}}

In this section, we recall the definition of a Manin triple as well as the link with the 
Lie bi-algebra structure. Then, we introduce the notion of twist and study their consequences.

Throughout this section, $\fp$ is a finite dimensional Lie algebra. Similar definitions can be given for
$\fp$ infinite but, in this case, the notion of duality must be tackled with care. 

\subsection{Lie bi-algebra and Manin triple. \label{sec:ma}}

In this subsection, 
we recall the definitions of a Lie bi-algebra and of a Manin triple
 (see \cite{CP} for a review). 
Let $\fg$ be a Lie algebra. A Lie bi-algebra structure on  
$\fg$ is a skew-symmetric linear map
$\delta_{\fg}:\fg\rightarrow\fg\otimes \fg$, the cocommutator, such that $\delta_{\fg}^*$ is a Lie bracket 
and $\delta_{\fg}$ is a 1-cocycle of $\fg$ with values in $\fg\otimes \fg$
\ben \label{cocy}
\delta_{\fg}([x,y])=x.\delta_{\fg}(y)-y.\delta_{\fg}(x)\;.
\een
The couple $(\fg,\delta_{\fg})$ is a Lie bi-algebra.

To get explicit examples of this structure, it is usually more convenient to deal with 
the notion of Manin triple.
A Manin triple is a triple of Lie algebras $(\fp,\fp_+,\fp_-)$ such that
\begin{enumerate}[(i)]
\item there is a non-degenerate symmetric bilinear form $\langle~,~\rangle_\fp$ on $\fp$,
invariant under the adjoint action of $\fp$;
\item $\fp_+$ and $\fp_-$ are Lie subalgebras of $\fp$;
\item $\fp=\fp_+ \oplus \fp_-$ as vector spaces;
\item $\fp_+$ and $\fp_-$ are isotropic for $\langle~,~\rangle_\fp$
(i.e.\ $\langle \fp_\pm,\fp_\pm\rangle_\fp=0$).
\end{enumerate}
Indeed, this definition allows one to get a Lie bi-algebra since there exists a one-to-one correspondence 
between the Manin triple $(\fp,\fp_+,\fp_-)$ and a
Lie bi-algebra structure for $\fp_+$. More precisely, the inner product $\langle~,~\rangle_\fp$ determines 
an isomorphism of vector spaces $\fp_+\simeq \fp_-^*$. Then, the Lie bi-algebra structure on $\fp_+$, i.e.\
the cocommutator $\delta_{\fp_+}:\fp_+\rightarrow \fp_+\otimes\fp_+$, is provided by 
the dual of the commutator $[~,~]_{\fp_-}:\fp_-\otimes\fp_-\rightarrow \fp_-$ (see e.g. \cite{CP} 
for a proof that $\delta_{\fp_+}=([~,~]_{\fp_-})^*$ is a 1-cocycle).
We give two examples of this correspondence in sections \ref{sec:Uq} and \ref{sec:Y}.

\subsection{Manin triple twists. \label{sec:mat}}

We give here one of the main new structures introduced in this paper.
\begin{definition} 
Let $(\fp,\fp_+,\fp_-)$ be a Manin triple. An invariant (resp. anti-invariant) Manin triple twist $\phi$ is an automorphism of $\fp$ 
satisfying
 \begin{enumerate}[(i)]
\item $\phi$ is an involution i.e. $\phi^2=Id$;
\item $\phi(\fp_\pm)=\fp_\pm$;
\item  $\langle \phi(x),y\rangle_\fp=+\langle x,\phi(y)\rangle_\fp$ 
(resp.  $\langle \phi(x),y\rangle_\fp=-\langle x,\phi(y)\rangle_\fp$), for any $x,y\in \fp$.
\end{enumerate}
\end{definition}
As a direct consequence of this definition, a Manin triple twist $\phi$ provides the symmetric space
decompositions 
\begin{equation}\label{eq:sys}
 \fp_\pm=\fk_\pm\oplus\fm_\pm\quad\text{with}\quad \phi(\fk_\pm)=+\fk_\pm
\quad\text{and}\quad \phi(\fm_\pm)=-\fm_\pm\;,
\end{equation}
and $\fp=\fk\oplus\fm$ with $\phi(\fk)=+\fk$ and $\phi(\fm)=-\fm$.
We recall that one gets
\begin{equation}\label{eq:in}
 [\fk,\fk]\subset \fk\;,\quad[\fk,\fm]\subset \fm\;\quad\text{and}\quad [\fm,\fm]\subset \fk\;.
\end{equation}
Similar relations hold for $\fk_\pm, \fm_\pm$.
In particular, the +1-eigenspaces $\fk$, $\fk_\pm$ are Lie algebras.
As explained previously (see section \ref{sec:ma}), a Manin triple $(\fp,\fp_+,\fp_-)$ allows one to
obtain a Lie bi-algebra structure for the Lie algebra $\fp_+$.
We are now interested if it is possible to get similar result for the Lie algebra 
$\fk_+$. 
We treat separately in the following subsections the cases when $\phi$ is an invariant or 
an anti-invariant Manin triple twist.

\subsubsection{Invariant Manin triple twist. \label{sec:it}}

Let $\phi$ be an invariant Manin triple twist. Using this invariance,
it is easy to show that
\begin{equation}\label{eq:is1}
 \langle \fk_+,\fm_-\rangle_\fp=0\quad\text{and}\quad
\langle \fm_+,\fk_-\rangle_\fp=0\;.
\end{equation}
We deduce that $\fk_+\simeq \fk_-^*$ as vector spaces. Finally, we can conclude that 
$(\fk,\fk_+,\fk_-)$ is a Manin triple with the inner product deduced from that of $\fp$. 
Therefore, $\fk_+$ is a Lie bi-algebra with the 
cocommutator $\delta_{\fk^+}$ dual of the Lie bracket $[~,~]_{\fk^-}$.
\begin{rmk} An important point to notice is that the 
cocommutator $\delta_{\fk_+}$, constructed by this way, is a map from $\fk_+$ to $\fk_+\otimes\fk_+$
and is different from the cocommutator $\delta_{\fp_+}$ (constructed from the Manin triple $(\fp,\fp_+,\fp_-)$) 
restricted to $\fk_+$ which is a map from
$\fk_+$ to $\fk_+\otimes\fk_+ + \fm_+\otimes\fm_+$.
\end{rmk}
An example of an invariant Manin triple twist is given in section \ref{sec:pb}.

\subsubsection{Anti-invariant Manin triple twist.\label{sec:AIMTT}}

Let $\phi$ be an anti-invariant Manin triple twist.
Because of the anti-invariance, we get, instead of (\ref{eq:is1}),
\begin{equation}\label{eq:is2}
 \langle \fk_+,\fk_-\rangle_\fp=0\quad\text{and}\quad
\langle \fm_+,\fm_-\rangle_\fp=0\;.
\end{equation}
It follows that $\fk_+\simeq \fm_-^*$ and $\fm_+\simeq \fk_-^*$ as vector spaces. 
Since $[\fm_-,\fm_-]\subset \fk_-$, $\fm^-$ is not a Lie algebra and 
we cannot deduce a Lie bi-algebra structure for $\fk_+$.
However, we can define a new type of structure.
\begin{definition}\label{def:tau}
Let $\phi$ be an anti-invariant Manin triple twist for $(\fp,\fp_+,\fp_-)$ which 
leads to the symmetric space decomposition (\ref{eq:sys}).
The linear map $\tau:\fk_+\rightarrow \fm_+ \otimes \fk_+$ is a left Lie bi-ideal 
structure for the couple $(\fk_+,\fm_+)$ if it is the dual of the
following action of $\fk_-$ on $\fm_-$
\begin{eqnarray}
 \tau^*~:~\fk_-\otimes \fm_- &\to&~~ \fm_-\\
 a\otimes y \quad&\mapsto&~~ [a,y]_{\fp_-}\;. \nonumber
\end{eqnarray}
\end{definition}
This definition makes sense since the map $\tau^*$ is well-defined 
knowing that $[\fk_-,\fm_-]_{\fp_-}\subset \fm_-$. 
We can define similarly a right Lie bi-ideal structure by the 
linear map $\tau':\fk_+\rightarrow \fk_+ \otimes \fm_+$ dual of the
action $\fm_-\otimes \fk_- \rightarrow \fm_-, y\otimes a \mapsto [y,a]_{\fp_-}$. 
This new structure we have introduced here is the main ingredient of this article.
We give various examples of it in sections \ref{sec:Ons}, \ref{sec:aOns}, \ref{sec:yp} and
\ref{sec:ym}.

\begin{rmk}
 Let us emphasize that the bi-ideal structures we have just introduced are closely related to 
the notion of coisotropic subalgebras \cite{Lu,Ci}. Indeed, the subalgebra $\fk_+$ with the map $\tau+\tau'$
is, by definition, a coisotropic subalgebra. Therefore, our construction provides a practical 
method to get some explicit examples of them. It would be interesting to link our method with the 
previous ones (see, for example, \cite{Za} and \cite{Oh} for their deformation).
\end{rmk}

\section{Deformation\label{sec:defo}}

In this section, we recall briefly the theory of the deformations for the 
enveloping algebras in order to get the quantized universal enveloping 
algebras (QUEA). We recall, in particular, the link between the Hopf structure of the
QUEA and the Lie bi-algebra structure. 
Then, we generalize this theory to be able to deform the Lie bi-ideal
structure. We show that, instead of Hopf algebra, we get 
coideal algebra.

\subsection{Quantized universal enveloping algebras.\label{sec:QUE}}  

In this section, we give a brief summary of the theory of the deformation
based on the review \cite{CP}.
Let $(\fg,\delta)$ be a Lie bi-algebra. 
Its universal enveloping algebra $\cU(\fg)$
may be equipped with a coproduct defined by 
$\Delta(x)=x\otimes 1 +1 \otimes x$ for $x\in\fg$ 
(the action of $\Delta$ is extended to $\cU(\fg)$ as an algebra homomorphism).
This coproduct allows one to extend the action of $\delta$
to  $\cU(\fg)$ by
\begin{equation}
 \delta(a_1a_2)=\delta(a_1)\Delta(a_2)+\Delta(a_1)\delta(a_2)\;.
\end{equation}
This extension makes $(U(\fg),\Delta,\delta)$ into a co-Poisson-Hopf structure.

By definition, the deformation of the Lie bi-algebra $(\fg,\delta)$ is the 
Hopf\footnote{To be concise, we do not mention the counit and the antipode but they 
can be introduced easily.} deformation $(\cU_\hbar(\fg),\Delta_\hbar)$ of 
its universal enveloping algebra $(\cU(\fg),\Delta,\delta)$
such that
\begin{equation}\label{eq:dD}
 \delta(x)=\frac{\Delta_\hbar(a)-\Delta^{\text{op}}_\hbar(a)}{\hbar}\quad (\text{mod } \hbar)
\end{equation}
 where $x\in \cU(\fg)$, $a\in \cU_\hbar(\fg)$ such that $x=a~~ (\text{mod }\hbar)$,
$\Delta^\text{op}_\hbar=\sigma\circ \Delta_\hbar$ and $\sigma$ flips both spaces.

The practical procedure, used in the following, to deform a Lie bi-algebra structure consists in
getting the simplest coproduct satisfying (\ref{eq:dD}) then to deform commutation relations
such that this coproduct is a homomorphism.
In sections \ref{sec:Uq} and \ref{sec:Y}, we give two well-known examples: 
the quantum universal enveloping algebra of the loop algebra based on $sl_2$ and the Yangian of $sl_2$.

\subsection{Coideal subalgebra. \label{sec:coia}} 

In this subsection, using the notations of section \ref{sec:MT},
we define the deformation
of the universal enveloping algebra $\cU(\fk_+)$ with 
the left Lie bi-ideal structure $\tau$.
Evidently, we cannot get a Hopf algebra but we obtain a coideal subalgebra.
Before going further, we give the definition of the left coideal subalgebra used in this paper.
A subalgebra $\cB$ of a Hopf algebra $(\cA,\Delta)$ is a left coideal subalgebra $(\cB,\fT)$ if
 $\fT$ is an algebra homomorphism from $\cB$ to $\cA\otimes\cB$, called coaction, 
with the coideal coassociativity property:
\begin{equation}\label{coass}
 (\Delta\otimes 1)\circ \fT=(1\otimes\fT)\circ \fT\;.
\end{equation}
We may define also a right coideal subalgebra.

As for Lie bi-algebra, we can equip the universal enveloping algebra $\cU(\fk_+)$
with a coproduct,
$\Delta(x)=x\otimes 1 +1 \otimes x$ for $x\in\fk_+$ 
which allows us to extend the action of $\tau$
to $\cU(\fk_+)$ by
\begin{equation}
 \tau(a_1a_2)=\tau(a_1)\Delta(a_2)+\Delta(a_1)\tau(a_2)\;,
\end{equation}
for $a_1,a_2\in\fk_+$.
\begin{rmk}
When we consider the universal enveloping algebra $\cU(\fk_+)$, the left bi-ideal 
structure $\tau$ maps $\cU(\fk_+)$ to  $\cU(\fp_+) \otimes \cU(\fk_+)$ and not 
to $\cU(\fm_+) \otimes \cU(\fk_+)$ in comparison with $\tau : \fk_+ \to \fm_+ \otimes \fk_+$.  
\end{rmk}
This extension does not 
make $(\cU(\fk_+),\Delta,\tau)$ into a co-Poisson-Hopf structure since
$\tau:\cU(\fk_+)\rightarrow \cU(\fp_+) \otimes \cU(\fk_+)$. Nevertheless,
the triple $(\cU(\fk_+),\Delta,\tau)$ shares properties with the co-Poisson-Hopf structure. 
But, even though we think it is an interesting problem, it is not in the 
scope of this article to give these properties.

We are now in a position to define the deformation for $\cU(\fk_+)$.
\begin{definition}
Let $(\fp_+,\delta)$ be the Lie bi-algebra defined in section \ref{sec:ma} and
$(\fk_+,\tau)$ be the Lie bi-ideal as in the definition \ref{def:tau}.
A deformation of $(\fp_+,\fk_+,\delta,\tau,\Delta)$ 
is $(\cU_\hbar(\fp_+),\cU_\hbar(\fk_+),\Delta_\hbar,\fT_\hbar)$ such that
\begin{enumerate}[(i)]
\item the Hopf algebra $(\cU_\hbar(\fp_+),\Delta_\hbar)$ is a deformation of the Lie bi-algebra 
$(\fp_+,\delta)$;
\item the algebra $\cU_\hbar(\fk_+)$ is a deformation of the algebra $\cU(\fk_+)$;
\item $(\cU_\hbar(\fk_+),\fT_\hbar)$ is a left coideal subalgebra of $(\cU_\hbar(\fp_+),\Delta_\hbar)$;
\item $\fT_\hbar(a)=\Delta(x)+\hbar \tau(x) ~~(\text{mod}~~\hbar^2)$ 
where $x\in \cU(\fk_+)$, $a\in \cU_\hbar(\fk_+)$ such that $x=a \text{ mod }(\hbar)$.
\end{enumerate}
\end{definition}
We can give a similar definition with $\tau'$, the right bi-ideal structure: we get a right 
coideal instead of a left one. It would be interesting to study the existence and the uniqueness of this type of 
deformations for any bi-ideal structures but it is already a complicated problem for the Lie bi-algebra structure.
Therefore, instead of giving a general proof, we provide in section \ref{sec:exM} 
different examples where such deformations exist.

\section{Examples of Manin triple twists and their deformation \label{sec:exM}} 

To illustrate the previous theoretical construction, we consider two different well-known 
examples of Manin triples: one related to the loop algebra of $sl_2$, denoted $\cL$, and 
another related to the half-loop algebra of $sl_2$, denoted $\cL^+$.
{}From these Manin triples, we obtain the Lie bi-algebra structures for $\cL$ and $\cL^+$ 
and we deform these algebras, following \cite{Dr,CP}.
For each of these Manin triples, we find different Manin triple twists, construct the
corresponding Lie bi-ideal structure and compute their deformation. These twists provide examples 
for the {\it invariant Manin triple twist} introduced in \ref{sec:it} and for the {\it anti-invariant 
Manin triple twist} introduced in \ref{sec:AIMTT}. 

We give for each case the defining relations in the ``Cartan type presentation'', 
in which it is easier to define the Manin triple (and the twist). 
Then, we obtain the Lie bi-algebra or the Lie bi-ideal structures. 
Finally, for simplicity, the deformation is given in the ``Serre type presentation''. 
The deformation in the ``Cartan type presentation'' is known for the 
Lie bi-algebras $\cL$ and $\cL^+$ (see e.g. \cite{Dr}). 
The Lie bi-ideal case will be considered elsewhere. 

\subsection{Deformation of the loop algebra ${\cL}$ and its twists. \label{sec:lie}}  

In this subsection, we deal with the loop algebra $\cL$ and its deformation: 
the quantum algebra $\cU_\hbar(\cL)$. Then, we consider two different twists of $\cL$.
The first one leads to the Onsager algebra $\cO$ \cite{Ons} and its deformation the q-Onsager algebra 
$\cO_\hbar$ \cite{B1,B2,IT1,IT2} 
whereas the second one leads to the augmented Onsager algebra $\overline{\cO}$ 
and its deformation the augmented q-Onsager algebra $\overline{\cO}_\hbar$ \cite{IT3}.  

\subsubsection{Quantum algebra $\cU_\hbar(\cL)$. \label{sec:Uq}}  

Here, for completeness, we recall the deformation of $\cL$ to get $\cU_\hbar(\cL)$. 
The loop algebra of $sl_2$, $\cL$, is the Lie algebra generated 
by $\{e_n,f_n,h_n|n\in \mathbb{Z}\}$ subject only to 
\ben\label{eq:CW}
&&[h_n,e_m]=2e_{n+m}\,, \quad [h_n,f_m]=-2f_{n+m}\,, \quad [e_n,f_m]=h_{n+m}\,\\
&& \mb{and} [h_n,h_m]=[e_n,e_m]=[f_n,f_m]=0\,, \nonumber
\een
 for $n,m\in \mathbb{Z}$. We use the inner product on $\cL$ defined by
\begin{equation}\label{innerL}
 \langle e_n,f_{-n} \rangle_\cL=-1\quad\text{and}\quad\langle h_n,h_{-n} \rangle_\cL=-2\;
\end{equation}
and vanishing otherwise.\\

Although the Lie bi-algebra structure for the loop algebra is well-known, we recall 
the main steps of its derivation from the Manin triple (see e.g. \cite{CP} chapter 1.3). 
We introduce $\cD=\cL\oplus\cL$ (as Lie algebras) and both subalgebras
\ben\label{MTD}
 \cD^+&=&\{(x,x)|x\in\cL\}\\
\label{MTD2}
 \cD^-&=&\{(e_n,0),(f_m,0),(0,f_{-n}),(0,e_{-m}),(h_n,0),(0,h_{-n}),(h_0,-h_0)|n> 0,m\geq 0\}\;.
\een
Then, $(\cD,\cD^+,\cD^-)$ is a Manin triple with the inner product 
\begin{equation}\label{pa}
 \langle (x,y), (x',y') \rangle_\cD=\langle x,x' \rangle_\cL -\langle y,y' \rangle_\cL\;.
\end{equation}
The subalgebra $\cD^+$ is isomorphic to $\cL$. Then, the Lie bi-algebra structure on $\cD^+$, i.e. 
the cocommutator $\delta : \cD^+ \to \cD^+ \otimes \cD^+$, 
dual of the commutator on the Lie algebra $\cD^-$ can be explicitly computed from the relation
\ben
\langle \delta\big((x,x)\big), (x',y')\otimes (x'',y'') \rangle_\cD= \langle (x,x), [(x',y'),(x'',y'')] \rangle_\cD
\een
for any $(x,x) \in \cD^+$ and $(x',y'), (x'',y'') \in \cD^-$.
As an example, we consider the case $(x,x)=(h_0,h_0)$. 
The non-vanishing term is given by $[(x',y'),(x'',y'')] =(h_0,-h_{0})$. 
This element cannot be obtained from commutation relations 
in $\cD^-$. It follows that $\delta\big((h_0,h_0)\big)=0$. 
With the same method, we obtain
\ben
&& \delta(h_0)=0\,,\quad
\delta(e_0)=\frac{1}{2}h_0\wedge e_0\,,\quad
\delta(f_0)=\frac{1}{2}h_0\wedge f_0\,,\quad
\delta(e_1)=-\frac{1}{2}h_0\wedge e_1\,,\quad 
\delta(f_{-1})=-\frac{1}{2}h_0\wedge f_{-1}. 
\een
To write the previous relations, we have used the isomorphism $\cD^+ \to \cL$ given by $(x,x)\mapsto x$
 and the wedge product $x\wedge y=x\otimes y-y\otimes x$.
Naturally, they satisfy cocycle condition (\ref{cocy}).\\

The knowledge of the cocommutator allows one to deform $\cU(\cL)$, 
as explained in section \ref{sec:QUE}. To simplify the computation of the deformation,
we use the ``Serre-Chevalley presentation'' of this algebra. 
It is given by the generators $\{X^\pm_i,H_i\;|\;i=0,1\}$ and the defining relations
\ben
&&[H_i,X^{\pm}_j]=\pm a_{ij}X^{\pm}_j\;, \quad [X^+_i,X^-_j]=\delta_{i,j}H_i\;,\quad
[X^{\pm}_i,[X^{\pm}_i,[X^{\pm}_i,X^{\pm}_{i\pm1}]]]=0\;, \quad H_0+H_1=0\;,
\een
with $A=\{a_{ij}\}$ the symmetric extended Cartan matrix of $\cL$ explicitly 
given by $ a_{ii}=2 $ and $ a_{ij}=-2$ for $i\neq j$.
The map from the second to the first formulation of $\cL$ is given by
\ben \label{sertocar}
H_1 \mapsto h _0\;,\quad H_0 \mapsto -h _0\;,\quad X^+_1 \mapsto e_0\;,\quad 
X^+_0 \mapsto f_{-1}\;,\quad X^-_1 \mapsto f_0\;,\quad X^-_0 \mapsto e_1\;.  
\een
Then, the Lie bi-algebra structure, in this presentation, is given 
by\footnote{The difference with the standard Lie bi-algebra structure 
is due to a factor $-\frac{1}{2}$ in the inner product (\ref{innerL}).}
\ben\label{eq:biLpUh}
 \delta(X^\pm_i)=\frac{1}{2}H_i\wedge X^\pm_i , \quad \delta(H_i)=0.
\een

The deformation of this algebra $\cU(\cL)$ from the previous cocommutator is considered in details in \cite{CP}. 
One obtains the quantum algebra $\cU_\hbar(\cL)$ with
the coproduct given by
\ben \label{DeltaUh}
\Delta_\hbar(\bH_i)=\bH_i\otimes 1+1 \otimes \bH_i,\quad
\Delta_\hbar(\bX^\pm_i)=\bX^\pm_i\otimes e^{-\frac{\hbar}{4}\bH_i}+e^{\frac{\hbar}{4}\bH_i} \otimes \bX^\pm_i,
\een
and the defining commutation relations 
\ben \label{comUh}
&&[\bH_i,\bX^\pm_j]=\pm a_{ij}\bX^\pm_j, \quad 
[\bX^+_i,\bX^-_j]=\delta_{i,j}\frac{e^{\frac{\hbar}{2}\bH_i}-e^{-\frac{\hbar}{2}\bH_i}}
{e^{\frac{\hbar}{2}}-e^{-\frac{\hbar}{2}}},\\
&&\null[\bX^\pm_i,[\bX^\pm_i,[\bX^\pm_i,\bX^\pm_{i\pm1}]_{(\hbar)}]_{(-\hbar)}]=0, \\
&& \bH_1+\bH_0=0,
\een
with $[x,y]_{(\hbar)}=xy-e^{\hbar}yx$. 
The bold symbols correspond to the deformations of the non-bold ones.


\subsubsection{q-Onsager algebra $\cO_\hbar$. \label{sec:Ons}} 
 
We consider now the subalgebra of the loop algebra $\cL$ invariant under the automorphism 
$\eta_1$ : $\cL \to \cL$ defined by
\begin{equation}\label{eq:eta1}
\eta_1(e_n)=f_{-n}~,~~\eta_1(f_n)=e_{-n}~~\text{ and }~~\eta_1(h_n)=-h_{-n}.
\end{equation}
This subalgebra is generated by $\{A_i,G_i\}$ with
\ben\label{eq:AG}
 A_i =2(e_i+f_{-i})~~\text{ and }~~G_i=h_i-h_{-i}. 
\een
It is easy to check that they satisfy the relations 
\ben \label{Ons}
[A_i,A_j]= 4 G_{i-j},\quad [A_i,G_j]=2(A_{i+j}-A_{i-j}), \quad
[G_i,G_j]=0, \quad G_i=-G_{-i}\;,
\een
which are the ones defining the Onsager algebra $\cO$ \cite{Ons} 
(see also \cite{Dav1,DaRo} for more details about the relation between $\cL$ and $\cO$). \\

To calculate the left bi-ideal structure for the Onsager algebra $\cO$, 
we need to extend the automorphism $\eta_1$ to $\cD=\cL\oplus\cL$ with the following map 
\begin{eqnarray}
 \phi_1:&\cD&\rightarrow~~ \cD \\
&(x,y)&\mapsto ~~(\eta_1(y),\eta_1(x)).\nonumber 
\end{eqnarray}
By direct computation, we show that $\phi_1$ is an anti-invariant Manin triple twist for the Manin triple 
$(\cD,\cD^+,\cD^-)$ defined by (\ref{MTD}), (\ref{MTD2}) and (\ref{pa}).
This twist implies 
the following symmetric space decompositions 
\ben
\cD^+=\ck^+\oplus \cm^+ &\text{with}
&\ck^+=\{(e_i+f_{-i},e_i+f_{-i}),(h_j-h_{-j},h_j-h_{-j})~|~i \in \mathbb{Z},j>0\}, \\
&&\cm^+= \{(e_i-f_{-i},e_i-f_{-i}),(h_k+h_{-k},h_k+h_{-k})~|~i \in \mathbb{Z},k\geq 0\}, \nonumber\\
\cD^-=\ck^-\oplus \cm^-&\text{with}
&\ck^-=\{(e_i,f_{-i}),(f_j,e_{-j}),(h_j,-h_{-j})~|~i>0,j\geq0\},\\
&&\cm^-=\{(e_i,-f_{-i}),(h_i,h_{-i}),(f_j,-e_{-j})~|~i>0,j\geq0\}\;.\nonumber
\een
We recover that 
$\langle \ck^+,\ck^-  \rangle_\cD= \langle \cm^+,\cm^-  \rangle_\cD=0 $ proven, in general, 
in section \ref{sec:AIMTT}.

We remark that the $+1$-eigenspace $\ck^+$ of $\cD^+$ is isomorphic to the Onsager 
algebra, with $(x,x)\mapsto x$. Then, we calculate the left Lie bi-ideal structure for it, 
following section \ref{sec:AIMTT}. This structure, i.e. $\tau : \ck^+ \to \cm^+ \otimes \ck^+$, 
dual of the commutator $[\;,\;] : \ck^- \otimes \cm^- \to  \cm^-$, can be explicitly computed from 
\ben
 \langle \tau \big((x,x)\big), (x',y')\otimes (x'',y'') \rangle_\cD= 
\langle (x,x), [(x',y'),(x'',y'')] \rangle_\cD
\een
with $(x,x)\in \ck^+$, $(x',y')\in \ck^-$ and $(x'',y'')\in \cm^-$.
After straightforward calculations, we obtain the left Lie bi-ideal structure for the Onsager algebra. 
The important values are
\ben\label{eq:tauu}
\tau(A_1)&=&-\frac{1}{2} h_0 \otimes A_1, \quad \tau(A_{0})=\frac{1}{2}h_0 \otimes A_{0}
\een
where we have used the map $(x,x) \mapsto x$ and relations (\ref{eq:AG}).\\

Similarly to the case of $\cL$, to simplify the deformation procedure, we use 
another presentation of the algebra $\cO$, called the Dolan-Grady algebra \cite{DG}. 
This presentation is given in terms of $\{A,A^*\}$ subject to the Dolan-Grady relations
\begin{equation}
 [A,[A,[A,A^*]]]=16 [A,A^*]\quad,\quad[A^*,[A^*,[A^*,A]]]=16[A^*,A]\;.
\end{equation}
The map with the first presentation is given by \cite{Dav2,DaRo}
\ben
A^*\mapsto A_0 \mb{and} A\mapsto A_1.
\een
Let us notice that, from relations (\ref{sertocar}) 
and (\ref{eq:AG}), we can identify $A$ with $2(X_0^+ + X_0^-)$ and $A^*$ with $2(X_1^+ + X_1^-)$.  
{}From relation (\ref{eq:tauu}) and using also the map (\ref{sertocar}), the left Lie bi-ideal 
structure, in this presentation, is given by
\ben
\tau(A)=\frac{1}{2}H_0 \otimes A\;, \quad
\tau(A^*)=\frac{1}{2}H_1 \otimes A^*. 
\een

{}From this Lie bi-ideal structure and section \ref{sec:defo}, it is natural 
to look for the left coaction $\fT_\hbar :  \cO_\hbar \to \cU_{\hbar}(\cL)\otimes  \cO_\hbar$ in the form
\ben
&&\fT_\hbar(\bA)=\widetilde{\bA}\otimes 1+f(\bH_0) \otimes \bA\\
&&\mb{with} \widetilde{\bA}=2(\bX^+_0+\bX^-_0)+O(\hbar) \mb{and} f(x)=1+\frac{\hbar}{2}x+O(\hbar^2)\label{expA}.
\een
Then, using the coideal coassociativity (\ref{coass}), we obtain that
\ben
\Delta_\hbar (\widetilde{\bA})&=&\widetilde{\bA}\otimes 1+f(H_0)\otimes \widetilde{\bA} \label{ansDA}\\
\Delta_\hbar(f(H_0))&=&f(H_0)\otimes f(H_0)\;.
\een
The last equation means that $f(H_0)$ is group-like. 
Then, together with the expansion of $f(x)$ (\ref{expA}), 
it implies that $f(x)=e^{\frac{\hbar}{2}x}$ (see \cite{CP} chapter 6.4).
Finally, using (\ref{expA}), (\ref{ansDA}) and the coproduct of $\cU_\hbar(\cL)$ (\ref{DeltaUh}),
we get $\widetilde{\bA}=\rho_\hbar(\bX^+_0+\bX^-_0)e^{\frac{\hbar}{4}\bH_0}$ with $\rho_\hbar=2+O(\hbar)$.
We choose  $\rho_\hbar=[2]_\hbar=\frac{e^{\hbar}-e^{-\hbar}}{e^{\frac{\hbar}{2}}-e^{-\frac{\hbar}{2}}}$. 
Similar computations for $\bA^*$ can be done and finally we obtain the coaction
\ben
\fT_\hbar(\bA)&=& [2]_\hbar(\bX^+_0+\bX^-_0)e^{\frac{\hbar}{4}\bH_0}\otimes 1
+e^{\frac{\hbar}{2}\bH_0}\otimes \bA\;,\\
\fT_\hbar(\bA^*)&=&[2]_\hbar(\bX^+_1+\bX^-_1)e^{\frac{\hbar}{4}\bH_1}\otimes 1
+e^{\frac{\hbar}{2}\bH_1}\otimes \bA^*\;.
\een
To require that this coaction is an homomorphism of algebra implies the following 
deformation for the Dolan-Grady relations \cite{B2}
\ben
[\bA,[\bA,[\bA,\bA^*]_{(\hbar)}]_{(-\hbar)}]&=&([2]_\hbar)^4\,[\bA,\bA^*]\,\\
\null [\bA^*,[\bA^*,[\bA^*,\bA]_{(\hbar)}]_{(-\hbar)}]&=&([2]_\hbar)^4\,[\bA^*,\bA],
\een
with $[x,y]_{(\hbar)}=xy-e^{\hbar}yx$. 
These relations have been previously introduced in \cite{IT1,IT2,B1} and the algebra 
generated by these relations is called q-Onsager algebra.

Note that the generators $\widetilde{\bA}=\rho_\hbar(\bX^+_0+\bX^-_0)e^{\frac{\hbar}{4}\bH_0}$
and ${\bA}=\rho_\hbar(\bX^+_1+\bX^-_1)e^{\frac{\hbar}{4}\bH_1}$ already appear in \cite{DM} 
and correspond to conserved non-local 
charges of the quantum affine Toda field theory on the half-line related to $\cU_q(\widehat{sl}_2)$.

\subsubsection{Augmented q-Onsager algebra $\overline{\cO}_\hbar$. \label{sec:aOns}} 

We consider the subalgebra of $\cL$ invariant under the automorphism $\eta_2 : \cL \to \cL$ defined by
\begin{equation}\label{eq:eta2}
\eta_2(e_n)=e_{-n+1}~,~~\eta_2(f_{n})=f_{-n-1}~~\text{ and }~~\eta_2(h_n)=h_{-n}.
\end{equation}
This subalgebra is generated by $\{B_i,\bar{B}_i,K_i\}$ with
\ben
B_i=e_i+e_{-i+1},\quad \bar{B}_i=f_i+f_{-i-1} \mb{and} K_i=h_i+h_{-i}.
\een
It easy to check that they satisfy the relations 
 \ben\label{augOns}
&&[B_i,B_j]=[\bar{B}_i,\bar{B}_j]=[K_i,K_j]=0, \quad B_i=B_{-i+1}, \quad \bar{B}_i=\bar{B}_{-i-1}, 
\quad  K_i=K_{-i}\\
&&\null [B_i,\bar{B}_j]=(K_{i+j}+K_{-i+j+1}), \quad
\null [K_i,\bar{B}_j]=-2(\bar{B}_{i+j}+\bar{B}_{-i+j}),\quad  [K_i,B_j]=2(B_{i+j}+B_{-i+j}).
\een
We call this subalgebra the augmented Onsager algebra $\overline{\cO}$ as to
reference to the q-deformed case introduced in \cite{IT3}; see below.
\\

As for the Onsager case, we introduce a twist for the Manin triple $(\cD,\cD^+,\cD^-)$ defined by (\ref{MTD}) and (\ref{MTD2}).
For that, we extend the automorphism $\eta_2$ to the following map 
\begin{eqnarray}
 \phi_2:&\cD&\rightarrow~~ \cD  \\
&(x,y)&\mapsto ~~(\eta_2(y),\eta_2(x)).\nonumber
\end{eqnarray}
We show that it is an anti-invariant Manin triple twist for $(\cD,\cD^+,\cD^-)$.
In this case, the symmetric space decompositions 
$\cD^+=\ck^+\oplus \cm^+$ and $\cD^-=\ck^-\oplus \cm^-$ are given by
\ben
\ck^+&=&\{(e_i+e_{-i+1},e_i+e_{-i+1}),(f_j+f_{-j-1},f_j+f_{-j-1}),(h_j+h_{-j},h_j+h_{-j})\} \\
\cm^+&=& \{(e_i-e_{-i+1},e_i-e_{-i+1}),(f_j-f_{-j-1},f_j-f_{-j-1}),(h_i-h_{-i},h_i-h_{-i})\}
\\
\ck^-&=&\{(f_j,f_{-j-1}),(e_i,e_{-i+1}),(h_i,h_{-i})\}
\\
\cm^-&=&\{(f_j,-f_{-j-1}),(e_i,-e_{-i+1}),(h_j,-h_{-j})\}
\een
with $i>0$ and $j\geq0$.
The $+1$-eigenspace of $\cD^+$ is isomorphic to the augmented Onsager algebra.
Using the same procedure as for the Onsager case, we can obtain the left Lie bi-ideal structure 
of the augmented Onsager algebra.
The important values are
\ben
\tau(B_1)=\frac{1}{4}(e_1-e_{0})\otimes K_0\,,\quad
\tau(\bar B_0)=-\frac{1}{4}(f_{0}-f_{-1})\otimes K_0\mb{and} \tau(K_0)=0,
\een
where we have used the map $(x,x) \mapsto x$.\\

To simplify the quantification, we used another presentation of this subalgebra 
generated by $B,B^*,K$ satisfying the relations
\ben\label{aonser}
&& [B,[B,[B,B^*]]]=0, \quad
  \null [B^*,[B^*,[B^*,B]]]=0,\quad [K,B^*]=-4\,B^*, \quad  [K,B]=4\,B.
\een
The map with the first presentation is given by
\ben
B \mapsto B_1\;, \quad B^* \mapsto \bar{B}_0\;, \mb{and} K \mapsto K_0\;.
\een
The left Lie bi-ideal structure in this presentation is given by (using also (\ref{sertocar}))
\ben
\tau(B)=\frac{1}{4}(X^-_0-X^+_1)\otimes K, \quad 
\tau(B^*)=-\frac{1}{4}(X^-_1-X^+_0)\otimes K \mb{and}
\tau(K)=0.
\een

{}From this left Lie bi-ideal structure, 
we show that the left coaction of $\overline{\cO}_\hbar$ can be chosen as 
\ben
&&\fT_\hbar(\bK)=(\bH_1-\bH_0)\otimes 1 +1\otimes \bK, \\
&&\fT_\hbar(\bB)= \bX^+_1 e^{\frac{\hbar}{4}\bH_0} \otimes  
e^{-\frac{\hbar}{4}\bK}+\bX^-_0 e^{\frac{\hbar}{4}\bH_1}  
\otimes  e^{\frac{\hbar}{4}\bK}+1 \otimes \bB, \\
&&\fT_\hbar(\bB^*)= \bX^-_1 e^{\frac{\hbar}{4}\bH_0} 
\otimes  e^{-\frac{\hbar}{4}\bK}+\bX^+_0 e^{\frac{\hbar}{4}\bH_1} 
\otimes  e^{\frac{\hbar}{4}\bK}+1 \otimes \bB^*.
\een
To require that this coaction is an homomorphism of algebra implies the following deformation of (\ref{aonser})
\ben
&& [\bB,[\bB,[\bB,\bB^*]_{(\hbar)}]_{(-\hbar)}]=-\rho_\hbar\bB(e^{\frac{\hbar}{2}\bK}-e^{-\frac{\hbar}{2}\bK})\bB, \\
&&\null  [\bB^*,[\bB^*,[\bB^*,\bB]_{(\hbar)}]_{(-\hbar)}]=\rho_\hbar\bB^*(e^{\frac{\hbar}{2}\bK}-e^{-\frac{\hbar}{2}\bK})\bB^*,\\
&&\null [\bK,\bB^*]=-4\,\bB^*, \quad  [\bK,\bB]=4\,\bB,
\een
with $[x,y]_{(\hbar)}=xy-e^{\hbar}yx$ and $\rho_\hbar=(e^{\hbar}-e^{-\hbar})(e^{\hbar}+e^{-\hbar}+1)$. 
These relations have been previously introduced in \cite{IT2} and the 
corresponding algebra is called augmented q-Onsager algebra.

\subsection{Deformation of the half-loop algebra  $\cL^+$ and its twists. \label{sec:HLo}}

In this subsection, we deal with the half-loop algebra $\cL^+$ and its deformation: 
the Yangian $\cY_\hbar(sl_2)$. Then, we give an example of invariant Manin triple twist 
leading to an Hopf algebra and two different coideal subalgebras. 
Since the construction of the Manin triple twists and their quantification are similar 
to the previous section, we do not give all the details. 

\subsubsection{Yangian $\cY_\hbar(sl_2)$.\label{sec:Y}}

Here, for completeness, we recall the deformation of the half-loop algebra $\cL^+$ to get 
the Yangian $\cY_\hbar(sl_2)$\footnote{The algebra called Yangian is in fact $\cY_1(sl_2)$ 
but it is well-known that $\cY_\hbar(sl_2)\simeq \cY_{\hbar'}(sl_2)$ for $\hbar,\hbar'\neq 0$. 
So, by abuse, we call Yangian
any $\cY_\hbar(sl_2)$ for $\hbar\neq 0$.}.
The half-loop algebra $\cL^+$ is the Lie algebra, with unit, generated 
by $\{e_n,f_n,h_n|n\in \mathbb{N}\}$ subject only to 
\ben\label{eq:CWH}
&&[h_n,e_m]=2e_{n+m}\,, \quad [h_n,f_m]=-2f_{n+m}\,, \quad [e_n,f_m]=h_{n+m}\,\nonumber\\
&& \mb{and} [h_n,h_m]=[e_n,e_m]=[f_n,f_m]=0\,, 
\een
 for $n,m\in \mathbb{N}$.
The half-loop algebra $\cL^+$ may be seen as the subalgebra of $\cL$ defined in section \ref{sec:Uq}.
The Lie bi-algebra structure for the half-loop algebra is well-known \cite{Dr} however we recall 
briefly its construction from a Manin triple.
If $\cL^-$ denotes the subalgebra of $\cL$ with the elements $\{e_n,f_n,h_n|n< 0\}$,
one can show that 
$(\cL,\cL^+,\cL^-)$ is a Manin triple with the inner product $\langle\langle~,~ \rangle\rangle_\cL$ defined by
\begin{equation}\label{innerLm}
 \langle\langle e_n,f_{-n-1} \rangle\rangle_\cL=-1\quad\text{and}\quad
\langle\langle h_n,h_{-n-1} \rangle\rangle_\cL=-2\;
\end{equation}
and equals $0$ otherwise. Let us emphasize that this inner product for $\cL$ is different 
from the one (\ref{innerL}) used to construct $\cU_\hbar(\cL)$.
Using this Manin triple and following section \ref{sec:ma}, one can obtain a Lie bi-algebra structure 
for $\cL^+$. The cocommutator for the first modes are given explicitly by 
\ben\label{eq:biLp1}
 \delta(x_0)=0~~,~\delta(e_1)=h_0\wedge e_0~,~\delta(f_1)=f_0\wedge h_0\quad\text{and}
\quad\delta(h_1)=2e_0\wedge f_0~,
\een
for $x_0\in\{e_0,f_0,h_0\}$.\\
  
As in the previous cases, it is more convenient to deal with the ``Serre-Chevalley presentation'' when we 
want to deform an algebra. It is known that $\cL^+$ may be generated only with $\{e,f,h,J(e),J(f),J(h)\}$
subject to relations 
\ben
&&[e,f]=h\;, \qquad [h,f]=-2f\;,\qquad [h,e]=2e \\
&&J([x,y])=[x,J(y)]\;,\qquad J(a x+b y)=aJ(x)+bJ(y)
\een
and the following Serre type relation 
\begin{equation}\label{ter1}
 [[J(e),J(f)],J(h)]=0\;.
\end{equation}
In this presentation, the whole set of
elements are obtained as imbricated commutators. 
The map between both presentations is given by
\ben
e \mapsto e_0,\quad f \mapsto f_0,\quad h \mapsto h_0,\quad J(e) \mapsto e_1,\quad 
J(f) \mapsto f_1,\quad J(h)  \mapsto h_1.
\een
The Lie bi-algebra structure in  the ``Serre-Chevalley presentation'' is easily deduced from the previous map 
and from relations (\ref{eq:biLp1}) and can be written compactly as follows
\ben\label{eq:biLp2}
 \delta(x)=0~~,~\delta(J(x))=[x\otimes 1,t],
\een
where $x\in\{e,f,h\}$ and $t=e\otimes f+f\otimes e+\frac{1}{2} h\otimes h$.\\

We can now consider the deformation of the Lie bi-algebra $(\cL^+,\delta)$ 
to get the Yangian of $sl_2$ \cite{Dr}.
The easiest way to get a Hopf structure satisfying relation (\ref{eq:dD}) with the Lie bi-algebra 
structure given by (\ref{eq:biLp2}) is to define the coproduct of $\cY_\hbar(sl_2)$ as follows
\ben
\Delta_\hbar\quad:\qquad x&\mapsto&~~ x\otimes 1+1\otimes x\\
\bJ(x)&\mapsto&~~ \bJ(x)\otimes 1+1\otimes \bJ(x)+\frac{\hbar}{2}[x\otimes 1,t]\nonumber
\een
In order to make $\Delta_\hbar$ a homomorphism, we must deform the commutation relations. 
It follows that $\cY_\hbar(sl_2)$ is the associative unital algebra
with six generators $e,f,h,\bJ(e),\bJ(f),\bJ(h)$ and the following defining relations
\ben
&&[e,f]=h\;,\qquad [h,e]=2e\;,\qquad [h,f]=-2f,\\
\null&&[x,\bJ(y)]=\bJ([x,y])\;,\qquad \bJ(a x+b y)=a\, \bJ(x)+b\, \bJ(y), \label{simpsl2}
\een
where $x,y\in\{e,f,h\}$, $a,b\in\CC$, and 
\ben
\big[[\bJ(e),\bJ(f)],\bJ(h)\big]~=~\hbar^2~ \big(\bJ(e)f-e\, \bJ(f)\big) h. \label{terrsl2}
\een
This last commutation relation is known as the ``terrific relation'' \cite{Dr} 
(see also \cite{Mo}).

\subsubsection{Positive Borel of $\cU_\hbar(A_1^{(2)})$.\label{sec:pb}}

We consider now the subalgebra $\cT_1$ of the half-loop algebra $\cL^+$ invariant under
the automorphism 
$\varphi_1:\cL^+\rightarrow\cL^+$ defined by 
\begin{equation}\label{eq:phi}
\varphi_1(e_n)=(-1)^{n+1} e_n~,~~\varphi_1(f_n)=(-1)^{n+1} f_n~~\text{ and }~~\varphi_1(h_n)=(-1)^{n}h_n.
\end{equation}
This subalgebra $\cT_1$ is generated by $\{h_{2n},e_{2n+1},f_{2n+1}|n\geq 0\}$ with the commutation 
relations (\ref{eq:CWH}) restricted to these generators. We extend the definition of $\varphi_1$ to $\cL$
by considering (\ref{eq:phi}) valid for $n\in \Z$. 

We want now to compute a Lie bi-algebra structure for $\cT_1$ using an invariant Manin triple twist.
To obtain that, we need to define the following Lie algebra 
$\cA=\cL\oplus \mathfrak{h}$ (direct sum of Lie algebras)
where $\mathfrak{h}$ is the abelian Lie algebra generated by $\{h_0\}$.
We introduce the following subalgebras of $\cA$   
\begin{eqnarray}
 \cA^+&=&\{(h_0,h_0),(e_n,0),(h_m,0),(f_m,0)|n\geq 0,m>0\}\\
 \cA^-&=&\{(h_0,-h_0),(e_n,0),(h_n,0),(f_m,0)|n< 0,m\leq 0\}\;.
\end{eqnarray}
Then, $(\cA,\cA^+,\cA^-)$ is a Manin triple with the inner product
\begin{equation}
 \langle (x,h),(x',h') \rangle_\cA =\langle x,x'\rangle_{\cL}-\langle h,h'\rangle_{\cL}\;.
\end{equation}
with $\langle ~,~\rangle_{\cL}$ defined in (\ref{innerL}).
Then, we can show that the map 
\begin{eqnarray}
 \phi_\cA\quad:\qquad\cA&\rightarrow& \cA\\
(x,h)&\mapsto& (\varphi_1(x),\varphi_1(h))\;\nonumber
\end{eqnarray}
is an invariant Manin triple twist for $(\cA,\cA^+,\cA^-)$.
{}From this map, we get a symmetric space decomposition for $\cA,\cA^+$ and $\cA^-$
where the $+1$-eigenspace of $\cA^+$ is isomorphic to $\cT_1$.
Explicitly the isomorphism $\cA^+\to \cT_1$ is given by $(h_0,h_0)\mapsto h_0$, 
$(e_1,0)\mapsto e_1$ and $(f_1,0)\mapsto f_1$.
As explained, in the section \ref{sec:it}, an invariant Manin triple twist provides a 
cocommutator for the $+1$-eigenspace. In this case, this construction gives the following 
cocommutator for $\cT_1$ 
\begin{equation}\label{eq:dinv}
 \delta(h_0)=0~,~~\delta(e_1)=\frac{1}{2}(h_0\otimes e_1 -e_1\otimes h_0)
~,~~\delta(f_1)=\frac{1}{2}(f_1\otimes h_0- h_0\otimes f_1)\;. 
\end{equation}
We give only the values of the cocommutator for $h_0,~e_1$ and $f_1$ since we need only these  
elements in the following.
\\

To study the deformation of this Lie bi-algebra $\cT_1$, we introduce, as for the previous examples,
the ``Serre-Chevalley presentation'' for $\cT_1$.
The Lie algebra $\cT_1$ can be equivalently defined only with $\{H,E,F\}$ subject to
\begin{eqnarray}
\label{eq:T1}
&& [H,E]=2E\quad,\quad[H,F]=-2F\\
\label{eq:T2}
&&[E,[E,[E,F]]]=0\quad,\quad[F,[F,[F,E]]]=0\;.
\end{eqnarray}
Let us emphasize that the previous terrific relations for $\cT_1$ are on the level 4 
in comparison to the one for $\cL^+$ (\ref{ter1}) which 
 is on the level 3. For the half-loop algebra based on the other Lie algebras as 
$sl_n$ ($n>2$), $so_n$ or $sp_n$, the terrific relations are on the level 2.
The map between both presentations is given by
\ben
H \mapsto h_0\,,\quad E \mapsto e_1 \mb{and} \quad F \mapsto f_1\;.
\een

We denote by $\bH$, $\bE$ and $\bF$ the deformed generators 
corresponding to $H$, $E$, $F$ in $\cT_1$.
Following the lines of the standard deformation of the Lie algebra $sl_2$ (see chapter 6.4.A in \cite{CP}),
we may choose the following coproduct
\ben
\Delta_\hbar( \bH)&=& \bH\otimes 1+1\otimes \bH,\\
\Delta_\hbar( \bE)&=& \bE\otimes e^{\hbar \bH/2}+1\otimes \bE,\\
\Delta_\hbar( \bF)&=& \bF\otimes 1+e^{\hbar \bH/2}\otimes \bF.
\een
It is easy to verify that this coproduct satisfies relation (\ref{eq:dD}) with the Lie bi-algebra structure defined by
(\ref{eq:dinv}) and induces the following commutation relations, deformations of relations (\ref{eq:T1})
and (\ref{eq:T2})
\begin{eqnarray}
&&[\bH,\bE]=2\bE\quad\text{,}\quad[\bH,\bF]=-2\bF\,,\\
&&[\bE,[\bE,[\bE,\bF]]_{(-\hbar)}]_{(-2\hbar)}=0\,,\quad [\bF,[\bF,[\bF,\bE]]_{(-\hbar)}]_{(-2\hbar)}=0.
\end{eqnarray}
with $[x,y]_{(\hbar)}=xy-e^{\hbar}yx$. 

To identify this Hopf algebra, let us recall that 
the twisted quantum affine groups ${\cU_\hbar(A_n^{(2)})}$ are well-studied
Hopf algebras \cite{J2} and their positive modes defined a subalgebra, called quantum positive Borel.
Usually, one restricts the study of these algebras to the cases $n\geq 2$ since one needs exterior 
automorphisms of $A_n$: 
the twist by an inner automorphism provides an algebra which is isomorphic to the one we start with
(see e.g. \cite{GO}).
However, this last statement is true only for the whole loop algebra and is no longer valid for the 
half-loop algebra. 
Therefore, by extension, we called the Hopf algebra we obtained the positive Borel
of ${\cU_\hbar(A_1^{(2)})}$.

\subsubsection{Coideal subalgebra $\cY^+_\hbar(sl_2)$ 
of the Yangian $\cY_\hbar(sl_2)$. \label{sec:yp}} 

In this subsection, we use again the subalgebra $\cT_1$ introduced in section \ref{sec:pb}
but this time we associate to it a Lie bi-ideal structure. To do that, we use the Manin triple
 $(\cL,\cL^+,\cL^-)$ with the inner product
$\langle\langle~,~\rangle\rangle_\cL$ (introduced in section \ref{sec:Y}) and we show that
the map $\varphi_1$, defined by (\ref{eq:phi}) (for $n\in \Z$), 
is an anti-invariant Manin triple twist.
Then, remarking that the $+1$-eigenspace of $\cL^+$ is isomorphic to $\cT_1$,
we can deduce a left bi-ideal structure for $\cT_1$. For example, we get
\begin{equation}\label{eq:toto}
 \tau(h_0)=0~,~~\tau(e_1)=-e_0 \otimes h_0~,~~\tau(f_1)=f_0 \otimes h_0 \;.
\end{equation}

\begin{rmk}
Let us emphasize that, although the algebras introduced in this example and in section \ref{sec:pb} 
are isomorphic, the additional structures associated are completely different. 
Indeed, for the former, we introduce the left Lie bi-ideal structure (\ref{eq:toto})
whereas, for the latter, we associate a Lie bi-algebra structure (\ref{eq:dinv}).
The difference in these structures leads to nonequivalent deformations.
\end{rmk}

To study the coideal subalgebra ${\cY^+_\hbar(sl_2)}$ of the Yangian $\cY_\hbar(sl_2)$,
deformation of the algebra $\cT_1$ with (\ref{eq:toto}), we use
the ``Serre-Chevalley presentation'' for $\cT_1$ given by (\ref{eq:T1}) and (\ref{eq:T2}).
We denote by $\bH$, $\bE$, $\bF$ the deformed generators in $\cY^+_\hbar(sl_2)$
corresponding, respectively, to the generators $H$, $E$, $F$ in $\cT_1$.
By the deformation procedure, we obtain the following left coaction 
\ben
\fT_\hbar(\bH)&=&h \otimes 1+1 \otimes \bH\\
\fT_\hbar(\bE)&=&\Big(\bJ(e)-\hbar~\frac{h\,e}{2}\Big) \otimes 1+ 1\otimes \bE-\hbar~ e \otimes \bH \\
\fT_\hbar(\bF)&=&\Big(\bJ(f)+\hbar~\frac{h\,f}{2}\Big) \otimes 1+1\otimes \bF+\hbar~ f \otimes \bH  
\een
and the following defining relations
\ben\label{eq:yp1}
[\bH,\bE]=2\,\bE,&& [\bH,\bF]=-2\,\bF,\\
\label{eq:yp2}
\Big[\bE,\big[\bE,[\bE,\bF]\big]\Big]= -12~\hbar^2~\,\bE\,\bH\,\bE,&& \Big[\bF,\big[\bF,[\bF,\bE]\big]\Big]=12~\hbar^2~\,\bF\,\bH\,\bF \label{TerYp}
\een
The embedding $\Psi_\hbar:\cY^{+}_\hbar(sl_2) \to \cY_\hbar(sl_2)$ is given by $\Psi_\hbar(\bH)=h$, 
$\Psi_\hbar(\bE)=\bJ(e)-\hbar\frac{he}{2}$ and $\Psi_\hbar(\bF)=\bJ(f)+\hbar\frac{hf}{2}$.\\

This coideal algebra $\cY^+_\hbar(sl_2)$ is isomorphic to the orthogonal twisted Yangian $\cR^+_\hbar(sl_2)$
studied previously in \cite{Ols}
where it has been defined via the so-called FRT presentation using the reflection equation. 
A detailed proof of this isomorphism is not explicitly given.
The proof is based on two principal steps. First, one proves that $\cY^+_\hbar(sl_2)$ 
is embedded in $\cR^+_\hbar(sl_2)$ by identifying elements of $\cR^+_\hbar(sl_2)$ 
satisfying the defining relations (\ref{eq:yp1}) and (\ref{eq:yp2}) of $\cY^+_\hbar(sl_2)$.
Second, by remarking that both algebras have a filtration, one can prove that 
both associated graded algebras are isomorphic (they are isomorphic to $\cT_1$).
These are sufficient to prove the isomorphism between $\cY^+_\hbar(sl_2)$ and $\cR^+_\hbar(sl_2)$.
A similar proof has been used in \cite{nG} for the twisted Yangians 
associated to $sl_n$ ($n\geq 3$).

To our knowledge, this presentation of the twisted Yangian $\cY^+_\hbar(sl_2)$ is new.
In fact, it was known that the three generators $\{h,~\bJ(e)-\hbar\frac{he}{2},~\bJ(f)+\hbar\frac{hf}{2}\}$ are
 the non-local charges of the Principal Chiral model on the half-line \cite{DMS} and
define a coideal subalgebra of $\cY(sl_2)$ but the terrific relations (\ref{TerYp}) were not known. 
Finally, let us recall that this coideal subalgebra is also isomorphic 
to the reflection algebra $\cB(2,1)$ \cite{MR}.

\subsubsection{Coideal subalgebra $\cY^-_\hbar(sl_2)$ 
of the Yangian $\cY_\hbar(sl_2)$. \label{sec:ym}} 

We consider now the subalgebra $\cT_2$ of the half-loop algebra $\cL^+$ invariant under
the automorphism 
$\varphi_2:\cL^+\rightarrow\cL^+$ defined by 
\begin{equation}\label{eq:chi}
\varphi_2(e_n)=(-1)^{n} e_n~,~~\varphi_2(f_n)=(-1)^{n} f_n~~\text{ and }~~\varphi_2(h_n)=(-1)^{n}h_n.
\end{equation}
This subalgebra $\cT_2$ is generated by $\{h_{2n},e_{2n},f_{2n}|n\geq 0\}$ and is isomorphic, as
Lie algebras, to $\cL^+$. We extend the definition of $\varphi_2$ to $\cL$
by considering (\ref{eq:chi}) valid for $n\in \Z$. 

The map $\varphi_2$ is an anti-invariant Manin triple twist for $(\cL,\cL^+,\cL^-)$ 
with the inner product
$\langle\langle~,~\rangle\rangle_\cL$.
We can show that, in this case, the $+1$-eigenspace of $\cL^+$ is $\cT_2$ 
and deduce a left Lie bi-ideal structure for $\cT_2$. 
For example, we get,
for $x_0\in\{e_0,f_0,h_0\}$,
\begin{equation}\label{eq:tau1}
 \tau(x_0)=0~~,~\tau(h_2)=2(e_1\otimes f_0-f_1\otimes e_0)~~,~
\tau(e_2)=h_1\otimes e_0-e_1\otimes h_0~~,~
\tau(f_2)=f_1\otimes h_0-h_1\otimes f_0\;.
\end{equation}

\begin{rmk}
This Lie bi-ideal structure (\ref{eq:tau1}) is different from the Lie bi-algebra structure 
defined by (\ref{eq:biLp1})
although the associated algebras in both cases are isomorphic. 
\end{rmk}

To study the coideal subalgebra ${\cY^-_\hbar(sl_2)}$ of the Yangian $\cY_\hbar(sl_2)$,
the deformation of the algebra $\cT_2$ with the structure studied previously, we use as usual
the ``Serre-Chevalley presentation'' for $\cT_2$. As noticed previously, $\cT_2$ is isomorphic to
$\cL^+$. Then $\cT_2$ may be generated only with $\{e,f,h,K(e),K(f),K(h)\}$
subject to relations 
\ben
&& [e,f]=h\;, \qquad [h,f]=-2f\;,\qquad [h,e]=2e\;,\\
&& K([x,y])=[x,K(y)]\;,\qquad K(a x+b y)=aK(x)+bK(y)
\een
and the following Serre type relation 
\begin{equation}\label{ter2}
 [[K(e),K(f)],K(h)]=0\;.
\end{equation}
The map between both presentations is given by
\ben
e \mapsto e_0,\quad f \mapsto f_0,\quad h \mapsto h_0,\quad K(e) \mapsto e_2,\quad 
K(f) \mapsto f_2,\quad K(h)  \mapsto h_2.
\een
and the left Lie bi-ideal structure if given by
\ben
 \tau(x)=0 \mb{and}  \tau(K(x))=[J(x)\otimes 1, t].
\een

We denote by $e$, $f$, $h$, $\bK(e)$, $\bK(f)$ and $\bK(h)$ the deformed generators in $\cY^-_\hbar(sl_2)$.
We follow the same procedure to deform $\cT_2$ than for the previous examples but, in this case, 
the results are more 
complicated since there exist non-vanishing 
terms in $\hbar^2$. 
We get the following left coaction $\fT_\hbar: \cY_\hbar^-(sl_2) \rightarrow \cY_\hbar(sl_2) \otimes \cY_\hbar^-(sl_2)$
defined by, for $x\in\{e,f,h\}$, 
\ben
\fT_\hbar(x)&=&x \otimes 1+1 \otimes x,\\
\fT_\hbar(\bK(h))&=&\Psi_\hbar(\bK(h))\otimes 1
+1\otimes \bK(h)  +\hbar\; [\bJ(h)\otimes 1, t]+
\frac{\hbar^2}{4}\Big([t,[t,h\otimes1]]+[[t,e\otimes 1],[t,f\otimes1]]\Big),\nonu
\fT_\hbar(\bK(e))&=&\Psi_\hbar(\bK(e))\otimes 1
+1\otimes \bK(e)  +\hbar\; [\bJ(e)\otimes 1, t]+
\frac{\hbar^2}{4}\Big([t,[t,e\otimes1]]+\frac{1}{2}[[t,h\otimes 1],[t,e\otimes1]]\Big),\nonu
\fT_\hbar(\bK(e))&=&\Psi_\hbar(\bK(e))\otimes 1
+1\otimes \bK(f)  +\hbar\; [\bJ(f)\otimes 1, t]+
\frac{\hbar^2}{4}\Big([t,[t,f\otimes1]]+\frac{1}{2}[[t,f\otimes 1],[t,h\otimes1]]\Big),\nonumber
\een
with $t=e\otimes f+f\otimes e+\frac{1}{2} h\otimes h$ and
\ben
\Psi_\hbar(
\bK(h))&=&
[\bJ(e),\bJ(f)]-\frac{\hbar}{4}[C,\bJ(h)]\\
\Psi_\hbar(
\bK(e))&=&
\frac{1}{2}[\bJ(h),\bJ(e)]-\frac{\hbar}{4}[C,\bJ(e)]
\\
\Psi_\hbar(
\bK(f))&=&\frac{1}{2}[\bJ(f),\bJ(h)]-\frac{\hbar}{4}[C,\bJ(f)]
\een
Then, the defining commutation relations of $\cY^-_\hbar(sl_2)$ are given by
\ben
&&[h,e]=2\,e\;,\qquad [h,f]=-2\,f\;,\qquad [e,f]=h\\
&&\bK([x,y])=[x,\bK(y)]\;,\qquad  \bK(a x+b y)=a \bK(x)+b \bK(y)
\een
with the terrific relation
\ben
\big[\bK(h),[\bK(e),\bK(f)]\big]&=&\hbar^2~\Big(\{\bK(f),\{\bK(h),e\}\}-\{\bK(e),\{\bK(h),f\}\}\Big)
+\frac{\hbar^4}{4} \,h\,[\bK(h),C]
\een
with $C=ef+fe+\frac{1}{2}h^2$ and $\{x,y\}=xy+yx$.

This coideal subalgebra is isomorphic to the symplectic twisted Yangian discovered in \cite{Ols} but
also to the reflection algebra $\cB(2,0)$ \cite{MR}. The proof is similar to the one for $\cY_\hbar^+(sl_2)$
(see section \ref{sec:yp}).

\section{Connection with Reflection algebra \label{sec:ReflEq}}

In this section, we used the previous coideal algebras to construct,
 from intertwiner technique \cite{J2,DMS,DM}, 
scalar $K$-matrices solutions of the reflection equation \cite{Cher,Sk}. 
It allows us to show the connection with known reflection algebras \cite{MRS,MR} 
characterized by the choice of the scalar $K$-matrices. 

\subsection{Scalar K-matrices for $\cO_\hbar$ and $\overline{\cO}_{\hbar}$.}

In this subsection, we construct scalar $K$-matrices for the coideal subalgebras of $\cU_\hbar(\cL)$: 
the q-Onsager algebra $\cO_\hbar$  (see section \ref{sec:Ons}) and
the augmented q-Onsager algebra $\overline{\cO}_{\hbar}$ (see section \ref{sec:aOns}).

The intertwiner equation for the scalar $K$-matrices is given by \cite{DM}, for any element $X$ of
 $\cO_\hbar$ ($\overline{\cO}_{\hbar}$),
\ben \label{intUh1}
K(u)(\pi_u \otimes \epsilon) \circ \fT_{\hbar}(X) = (\pi_{u^{-1}} \otimes \epsilon) \circ \fT_{\hbar}(X)K(u)
\een
with $\pi_u$ a representation of the quantum algebra $\cU_\hbar(\cL)$, depending on the spectral parameter $u$,
 and $\epsilon$ the counit of $\cO_\hbar$ ($\overline{\cO}_{\hbar}$). For $\cO_\hbar$, the counit is given by
\ben\label{eq:cuA}
\epsilon(\bA)=a, \quad \epsilon(\bA^*)=a^*
\een
with $a,a^*$ two arbitrary parameters and for $\overline{\cO}_\hbar$, we have 
\ben
\epsilon(\bB)=\epsilon(\bB^*)=0, \quad \epsilon(\bK)=c
\een
with $c$ an arbitrary parameter. \\

For example, we consider the case where $\pi_u$ is the fundamental evaluation representation 
of the quantum algebra $\cU_\hbar(\cL)$ in the principal picture given by 
\ben
\pi_u(X^+_1)&=& \alpha E_{21}\,, \quad \pi_u(X^-_1)=  \alpha^{-1} E_{12}\,, \quad \pi_u(H_1)=E_{11}-E_{22}\,, \\
\pi_u(X^+_0)&=&\alpha u E_{12}\,, \quad \pi_u(X^-_0)= \alpha^{-1} u^{-1}E_{21}\,, \quad \pi_u(H_0)=-E_{11}+E_{22}, 
\een
with $E_{ij}$ the two by two matrix with $1$ in the entry $(i,j)$ and $0$ otherwise. 
The parameter $\alpha$, entering in this representation, will be fixed by the intertwiner equation. 
After straightforward computations, we show that the solution $K(u)$ of the intertwiner 
equation is given by, for the q-Onsager algebra $\cO_\hbar$
\ben\label{Knd}
K(u) \propto \begin{pmatrix} 
a^* u+a &[2]_\hbar \frac{ (u-u^{-1})}{e^{\frac{\hbar}{2}}-e^{-\frac{\hbar}{2}}} \\
[2]_\hbar \frac{ (u-u^{-1})}{e^{\frac{\hbar}{2}}-e^{-\frac{\hbar}{2}}} & a^* u^{-1}+a
\end{pmatrix}
\een
with $\alpha=e^{\frac{\hbar}{4}}$ and by, for the augmented q-Onsager algebra $\overline{\cO}_\hbar$,
\ben\label{Kd}
K(u) \propto \begin{pmatrix} 
u+e^{\frac{(2 c-1) \hbar}{4}} & 0 \\
0 & u^{-1}+e^{\frac{(2 c-1) \hbar}{4}} 
\end{pmatrix}
\een
with $\alpha=e^{-\frac{\hbar}{8}}$. 
Both solutions satisfy the reflection equation 
\ben\label{REU}
R_{21}(u/v)K_1(u)R_{12}(uv)K_2(v)=K_2(v)R_{21}(uv)K_1(u)R_{12}(u/v)
\een
with $R_{12}(u)$ the trigonometric solution of the Yang Baxter equation \cite{J1} given by
\ben
R_{12}(u)= \begin{pmatrix} 
1& 0 & 0 &0 \\
0& e^{\frac{\hbar}{2}} \frac{1-u}{1-e^{\hbar}u} &  \frac{1-e^{\hbar}}{1-e^{\hbar}u}  & 0\\
0&u \frac{1-e^{\hbar}}{1-e^{\hbar}u} & e^{\frac{\hbar}{2}} \frac{1-u}{1-e^{\hbar}u}& 0 \\
0 & 0 & 0 & 1
\end{pmatrix}.
\een
Up to conjugation by an arbitrary invertible diagonal 
matrix and multiplication by a scalar (that leave the reflection equation invariant) the solution (\ref{Knd}), 
that also appears in \cite{DM}, corresponds to the one of \cite{GZ,dVGR}.

The diagonal solution (\ref{Kd}) may be obtained from multiplication of (\ref{Knd}) by a scalar $d$ 
and sending this scalar to zero keeping $d a^*=1$ and $da = e^{\frac{(2 c-1) \hbar}{4}}$.
Let us remark that this scaling limit is singular for the counit of the Onsager algebra 
(see equation (\ref{eq:cuA})). Therefore, to use the intertwiner relation to determine the diagonal solution,
we need the augmented q-Onsager algebra.
  
From another point of view, if one constructs reflection algebra using 
dressing procedure from a scalar $K$-matrix  
(i.e. with monodromy matrix of $\cU_{\hbar}(\cL)$ in the RLL presentation \cite{FST}) \cite{Sk,B1,B2,MRS}, 
the algebra obtained from the non-diagonal $K$-matrix (\ref{Knd}) corresponds
 to the q-Onsager algebra and from the diagonal 
$K$-matrix (\ref{Kd}) corresponds to the augmented q-Onsager algebra.

  
\subsection{Scalar K-matrices for $\cY^{\pm}(sl_2)$.}

In this subsection, we construct scalar $K$-matrices for the coideal subalgebras of $\cY_\hbar(sl_2)$: 
the twisted Yangians $\cY^+_\hbar(sl_2)$  (see section \ref{sec:yp}) and
$\cY^-_\hbar(sl_2)$ (see section \ref{sec:ym}).

As previously, the intertwiner equation for the scalar $K$-matrices is given by \cite{DM}, for any element $X$ of
 $\cY^+_\hbar(sl_2)$ ($\cY^-_\hbar(sl_2)$),
\ben \label{intUh}
K(u)(\pi_u \otimes \epsilon) \circ \fT_{\hbar}(X) = (\pi_{-u} \otimes \epsilon) \circ \fT_{\hbar}(X)K(u)
\een
with $\pi_u$ a representation of the quantum algebra $\cY_\hbar(sl_2)$, depending on the spectral parameter $u$,
 and $\epsilon$ the counit of $\cY^+_\hbar(sl_2)$ ($\cY^-_\hbar(sl_2)$). For $\cY^+(sl(2))$, the counit is
\ben
\epsilon(\bH)=c, \quad \epsilon(\bE)=\epsilon(\bF)=0 
\een
and, for $\cY^-(sl(2))$, it is
\ben
\epsilon(x)=\epsilon(\bK(x))=0 \mb{for} x \in \{e,f,h\}.
\een

For example, we consider the case where $\pi_u$ is the fundamental evaluation representation 
of the Yangian $\cY_\hbar(sl(2))$ given by 
\ben
\pi_u(e)=E_{21}\,, \quad \pi_u(f)=E_{12}\,, \quad \pi_u(h)=E_{11}-E_{22} \mb{and} \pi_u(\bJ(x))=(u+\alpha) \pi_u(x)
\een
with $\alpha$ an arbitrary parameter. 
After straightforward computations, we show that the solution $K(u)$ of the intertwiner equation is given by, 
for $\cY^+(sl(2))$ 
\ben \label{KsYp}
K(u) \propto \begin{pmatrix} 
1+\frac{c \hbar}{u} & 0 \\
0 & -1+\frac{c \hbar}{u}
\end{pmatrix}
\een
with $\alpha=1$ and by, for $\cY^-(sl(2))$,
\ben \label{KsYm}
K(u)\propto \begin{pmatrix} 
1& 0 \\
0 & 1
\end{pmatrix}
\een
with $\alpha=0$. 
Both solutions satisfy the reflection equation
\ben \label{REY}
R_{12}(u-v)K_1(u)R_{12}(u+v)K_2(v)=K_2(v)R_{12}(u+v)K_1(u)R_{12}(u-v)
\een
with $R_{12}(u)$ the rational solution of the Yang-Baxter equation given by
\ben
R_{12}(u)= \begin{pmatrix} 
1& 0 & 0 &0 \\
0&\frac{ u}{u+\hbar} & \frac{  \hbar}{u+\hbar}  & 0\\
0& \frac{  \hbar}{u+\hbar}  & \frac{ u}{u+\hbar}  & 0 \\
0 & 0 & 0 & 1
\end{pmatrix}\;.
\een
The general solution of the reflection equation (\ref{REY}) has been obtained by direct resolution 
in \cite{dVGR} 
and contains four arbitrary parameters. It also appears in \cite{DMS} from intertwiner technique. 
This solution is equivalent to our solution (\ref{KsYp}) related to $\cY^+(sl(2))$. 
The three supplementary parameters come from the invariance 
of (\ref{REY}) by the conjugation with an arbitrary invertible matrix
 $M\in End(\CC^2)$: $M_1M_2R_{12}(u)M^{-1}_1M^{-1}_2=R_{12}(u)$. 
Note that the identity solution (\ref{KsYm}), that also appears in \cite{DMS}, 
cannot be recovered from the solution of \cite{dVGR}. 
It is consistent with the existence of two algebras inside the reflection algebra constructed from (\ref{REY}).  

The $K$-matrix obtained for $\cY^+(sl(2))$ corresponds to the one used to define $\cB(2,1)$ \cite{MR} 
or the orthogonal twisted Yangian \cite{Ols} (up to some conjugation and transposition) and 
the $K$-matrix obtained for $\cY^-(sl(2))$ corresponds to the one used to define $\cB(2,0)$ \cite{MR} 
or the symplectic twisted Yangian \cite{Ols}.

\section{Conclusion\label{sec:conc}}

In this paper, we define a new structure associated to Lie algebras, called Lie bi-ideal structures, 
and provide a way to obtain examples thanks to particular automorphisms of Manin triples. 
We show that this structure is the good framework to obtain coideal algebras provided these 
Lie bi-algebras be deformed.

There remain a lot of interesting open questions.
Lie bi-algebra structures are known to guarantee 
the existence and the uniqueness of the Hopf algebra. The proof is based on the Hopf algebra cohomology.
Therefore, it is certainly of great interest to develop a coideal algebra cohomology and 
study the implications for the bi-ideal structures introduced in this paper. 
{}From another point of view, the construction of an analog of the coboundary for the bi-ideal structure 
and the computation of classical $K$-matrices are interesting questions. 

Regarding applications to integrable models with boundaries, this presentation 
of coideal algebras is interesting for the construction of dynamical $K$-matrices 
(i.e. with algebraic entries) that are solutions of the reflection equations \cite{BK,BF}. 
They allow one to construct integrable models with degrees of freedom on the boundaries \cite{BK}. 
On the other hand, we expect that they can be used to construct the algebraic $Q$ Baxter operators for open models. 
For periodic models this construction is known \cite{BLZ}. 
These algebraic $Q$ Baxter operators appear, for example, in the study of the correlation functions \cite{JMS}.
\vspace{2mm}\\ 
{Acknowledgement:} The authors thank P.~Baseilhac, A.I.~Molev, E.~Ragoucy and P.~Roche for 
their careful reading of the manuscript.

\end{document}